\documentclass[11pt]{article}
\usepackage{geometry}
\usepackage[parfill]{parskip} 
\usepackage{graphicx}
\usepackage{amssymb}
\usepackage{epstopdf}
\usepackage{amsmath,amsthm,amscd,amssymb}
\usepackage{latexsym}
\usepackage[colorlinks,citecolor=red,pagebackref,hypertexnames=false]{hyperref}
\numberwithin{equation}{section}
\newtheorem{theorem}{Theorem}[section]
\newtheorem{corollary}[theorem]{Corollary}

\newtheorem{lemma}[theorem]{Lemma}
\newtheorem{proposition}[theorem]{Proposition}
\newtheorem{remark}[theorem]{Remark}
\newtheorem{problem}[theorem]{Problem}

\newtheorem{definition}[theorem]{Definition}

\title{Results on the Erd\H os-Falconer \break distance  problem   in $\mathbb{Z}_q^d$ for odd $q$}
\author{{\sc David J. Covert} \\
University of Missouri - St. Louis  \\
covertdj@umsl.edu}

\date{}

\begin{document}
\maketitle

\abstract{The Erd\H os-Falconer distance problem in $\mathbb{Z}_q^d$ asks one to show that if $E \subset \mathbb{Z}_q^d$ is of sufficiently large cardinality, then $\Delta(E) := \{(x_1 - y_1)^2 + \dots + (x_d - y_d)^2 : x, y \in E\}$ satisfies $\Delta(E) = \mathbb{Z}_q$.  Here, $\mathbb{Z}_q$ is the set of integers modulo $q$, and $\mathbb{Z}_q^d = \mathbb{Z}_q \times  \dots \times \mathbb{Z}_q$ is the free module of rank $d$ over $\mathbb{Z}_q$.  We extend known results in two directions.  Previous results were known only in the setting $q = p^{\ell}$, where $p$ is an odd prime, and as such only showed that all units were obtained in the distance set.  We remove the constriction that $q$ is a power of a prime, and despite this, shows that the distance set of $E$ contains \emph{all} of $\mathbb{Z}_q$ whenever $E$ is of sufficiently large cardinality.}

\section{Background and Results}

A large portion of geometric combinatorics asks one to show that if a set is sufficiently large, then it exhibits some specific type of geometric structure.  One of the best known such results is the Erd\H os-distance problem.  Let $f(n)$ be the minimum number of Euclidean distances determined by any set of $n$ points in $\mathbb{R}^d$.  The classical Erd\H os-distance problem asks one to show that there exists a constant $c$ such that
\[
f(n) \geq \left\{ 
\begin{array}{cl}
cn^{1 - \underline{o}(1)} & d = 2
\\
cn^{2/d} & d \geq 3
\end{array}
\right.
\]
As usual, we write $X = \underline{o}(Y)$ if $X/Y$ tends to zero as some common parameter of $X$ and $Y$ tends to infinity.  The case $d = 2$ was recently resolved by Guth-Katz (\cite{GK}) with an ingenious application of the polynomial method.  The conjecture remains open for $d \geq 3$.  See \cite{Erdos, GIS, KT, SV, TV} and the references therein for more background and a  thorough treatise of the problem.

A continuous analog of the Erd\H os-distance problem is due to Falconer (\cite{Mattila}).  For a compact set $E \subset [0,1]^d$, let $\Delta(E) = \{|x - y| : x, y \in E\}$ be the set of pairwise Euclidean distances in $E$.  Falconer showed (\cite{Falconer}) that if $\dim_H(E) > \frac{d+1}{2}$, then $\mathcal{L}^1(\Delta(E)) > 0$, where $\dim_H(\cdot)$ denotes Hausdorff dimension and $\mathcal{L}^1(\cdot)$ denotes $1$-dimensional Lebesgue measure.  Further, he constructed a compact set $E \subset [0,1]^d$ with $\dim_H(E) = \frac{d}{2}$ such that $\Delta(E)$ did not have positive Lebesgue measure.  This led him to conjecture that if $E \subset [0,1]^d$ is compact and $\dim_H(E) > \frac{d}{2}$, then $\mathcal{L}^1(\Delta(E)) > 0$.  The best know results in this direction are due to Wolff (in the plane) and Erdo\u{g}an (for $d \geq 3$) who showed that $\mathcal{L}^1(\Delta(E)) > 0$ whenever $E \subset [0,1]^d$ is compact with $\dim_H(E) > \frac{d}{2} + \frac{1}{3}$ (\cite{Erdogan, Wolff}).

A finite field analogue of the distance problem was first considered by Bourgain-Katz-Tao (\cite{BKT}).  Let $\mathbb{F}_q^d$ denote the $d$-dimensional vector space over the finite field with $q$ elements.  For $E \subset \mathbb{F}_q^d$, define
\[
\Delta(E) = \{\| x - y\| : x, y \in E\} \subset \mathbb{F}_q
\]
where 
\[
\| x - y \| = (x - y) (x - y)^t = (x_1 - y_1)^2 + \dots + (x_d - y_d)^2.
\]
Clearly, $\| \cdot \|$ is not a norm, though this notion of distance is preserved under orthogonal transformations.  More precisely, let $O_d(\mathbb{F}_q)$ denote the set of $d \times d$ orthogonal matrices with entries in $\mathbb{F}_q$.  One can readily check that for $O \in O_d(\mathbb{F}_q)$, we have $\| x \| = \| Ox \|$.  Once a suitable notion of distance has been established in $\mathbb{F}_q^d$, the problem proceeds just as before.
\begin{problem}[Erd\H os-Falconer Distance Problem]\label{alpha}
Find the minimal exponent $\alpha$ such that if $q$ is odd then there exists a constant $C$ (independent of $q$) so that for any set $E \subset \mathbb{F}_q^d$ of cardinality $|E| \geq C q^{\alpha}$, we have $\Delta(E) = \mathbb{F}_q$.
\end{problem}
Note that when $q = p^2$, then $\mathbb{F}_{p^2}$ contains a subfield isomorphic to $\mathbb{F}_p$, and hence,  there exists a set $E \subset \mathbb{F}_q^d$ such that $|E| = q^{d/2}$ and $|\Delta(E)| = \sqrt{q}$.  This shows that the exponent $\alpha$ from Problem \ref{alpha} cannot be less than $\frac{d}{2}$ in general.

\begin{remark}
The methods used to attack the Erd\H os-Falconer distance problem are of a much different flavor when $q$ is even.  For example, in $\mathbb{F}_2^d$, we can construct a large set that contains only $0$ in its distance set as follows.  Let $E \subset \mathbb{F}_2^d$ be the set of vectors which have an even number of nonzero components.  Then, $\|x - y \| = 0$ in $\mathbb{F}_2$ since
\[
\| x - y \| = (x_1 - y_1)^2 + \dots (x_d - y_d)^2 = x_1 + y_1 + \dots + x_d + y_d .
\]
Note also that $E$ has cardinality
\[
|E| = \sum_{k = 0}^{\left\lfloor \frac{d}{2} \right\rfloor} {d \choose 2k} = 2^{d-1}.
\]
Thus, we have explicitly constructed a set of size $|E| = 2^{d-1}$ such that $\Delta(E) = \{0\}$.  Then, taking any set of size $|E| > 2^{d-1}$, gives $\Delta(E) = \mathbb{F}_2$ by the pigeonhole principle.  This gives the sharp exponent in $\mathbb{F}_2^d$ in the strongest possible sense.

\end{remark}
We shall henceforth assume $q$ is odd.  Iosevich-Rudnev (\cite{IR}) gave the first explicit exponent for the Erd\H os-Falconer distance problem in $\mathbb{F}_q^d$:
\begin{theorem}\label{IR}
There exists a constant $C$ so that if $E \subset \mathbb{F}_q^d$ has cardinality $|E| \geq C q^{\frac{d+1}{2}}$ then $\Delta(E) = \mathbb{F}_q$.
\end{theorem}
It would be reasonable to expect that whenever $E \subset \mathbb{F}_q^d$ with $|E| \geq Cq^{\frac{d}{2}}$ for a sufficiently large constant $C$, then $\Delta(E) = \mathbb{F}_q$, in line with the Falconer distance problem.  However, it was shown in \cite{HIKR} that Theorem \ref{IR} is sharp in odd dimensions in the sense that the exponent $\frac{d+1}{2}$ cannot be replaced by any smaller value.  It may still be the case that $\frac{d}{2}$ is the proper exponent in even dimensions.  The only known improvement occurs in the case $d=2$, where it has been shown (\cite{BHIPR, CEHIK}) that if $E \subset \mathbb{F}_q^2$, then there exists a constant $C$ such that whenever $|E| \geq Cq^{4/3}$, then $|\Delta(E)| \geq cq$ for some $ 0 < c \leq 1$.  Note that the exponent $\alpha = 4/3$ is in line with Wolff's exponent for the Falconer distance problem.  See \cite{GIS, KS} and the references contained therein for more on the Erd\H os-Falconer problem and related results.

Despite the Erd\H os-distance problem having been resolved in $\mathbb{R}^2$, the Falconer distance problem is open in all dimensions, and the finite field analogue is open in all even dimensions.  To try and obtain a better understanding of why this is the case, the author along with Iosevich and Pakianathan extended (\cite{CIP}) the Erd\H os-Falconer distance problem to $\mathbb{Z}_q$, the integers modulo $q$.  Here we let $\mathbb{Z}_q ^d = \mathbb{Z}_q \times \dots \times \mathbb{Z}_q$ denote the free module of rank $d$ over $\mathbb{Z}_q$.  For $E \subset \mathbb{Z}_q^d$, define $\Delta(E) = \{ \| x - y \| : x, y \in E\} \subset \mathbb{Z}_q$, where as before $\|x - y \| = (x_1 - y_1)^2 + \dots + (x_d - y_d)^2$.  We obtain the following results in this setting.

\begin{theorem}\label{CIP1}
Suppose that $E \subset \mathbb{Z}_q^d$, where $q = p^{\ell}$ is a power of an odd prime.  Then there exists a constant $C$ such that $\Delta(E) \supset \mathbb{Z}_q^{\times} \cup \{ 0 \}$ whenever $|E| \geq C \ell(\ell + 1) q^{ \frac{(2 \ell - 1)d + 1}{2 \ell}}$.
\end{theorem}

This result is a nice extension of Theorem \ref{IR} in the sense that when $\ell = 1$, $\mathbb{Z}_{p^{\ell}}$ is a field, and the exponents match those of Theorem \ref{IR} exactly.  Since Theorem \ref{IR} is sharp in odd dimensions, then Theorem \ref{CIP1} is sharp in odd dimensions as well, at least in the case $\ell = 1$.  In \cite{CIP} it was shown that Theorem \ref{CIP1} is close to optimal in the sense that there exists a value $b = b(p)$ such that $|E| = b q^{\left( \frac{2 \ell - 1}{2 \ell} \right) d}$, and yet $\Delta(E) \cap \mathbb{Z}_q^{\times} = \emptyset$.  This shows that for these constructed sets $E$, we have$|\Delta(E)| \leq p^{\ell - 1} = \underline{o}(q)$.  

It is of interest to extend Theorem \ref{CIP1} to non-units in $\mathbb{Z}_q$ and to the case $q \neq p^{\ell}$.  This is the purpose of the article, and our main result is  the following.

\begin{theorem}\label{main}
Suppose that $q$ has the prime decomposition $q  = p_1^{\alpha_1} \dots p_k^{\alpha_k}$, where $2 < p_1 < \dots < p_k$ and $\alpha_i > 0$ for each $i = 1, \dots , k$.  Suppose that $E \subset \mathbb{Z}_q^d$ for some $d > 2$.  Let $\tau(q) = \sum_{d \mid q} 1$ be the number of positive divisors of $q$.  Then, there exists a constant $C$ such that $\Delta(E) = \mathbb{Z}_q$ whenever 
\[
|E| \geq C\tau(q) q^d p_1^{- \frac{d-2}{2}}.
\]
\end{theorem}
\begin{remark}
As $\tau(q) = \underline{o}(q^{\epsilon})$ for all $\epsilon > 0$, our result is always nontrivial for $d \geq 3$.\end{remark}
Theorem \ref{main} immediately implies the following.
\begin{corollary}
Suppose that $E \subset \mathbb{Z}_q^d$, where $q = p^{\ell}$ is odd and $d > 2$.  Then there exists a constant $C$ such that $\Delta(E) = \mathbb{Z}_q$ whenever $|E| \geq C (\ell + 1) q^{\frac{( 2 \ell - 1)d + 2}{2 \ell}}$.
\end{corollary}

\subsubsection{Fourier Analysis in $\mathbb{Z}_q^d$}

For $f : \mathbb{Z}_q^d \to \mathbb{C}$, we define the (normalized) Fourier transform of $f$ as
\[
\widehat{f}(m) = q^{-d}\sum_{x \in \mathbb{Z}_q^d} f(x) \chi(-x \cdot m)
\]
where $\chi(x) = \exp(2 \pi i x/q)$.  Since $\chi$ is a character on the additive group $\mathbb{Z}_q$, we have the following orthogonality property.
\begin{lemma}\label{orthog} We have
\[
q^{-d} \sum_{x \in \mathbb{Z}_q^d} \chi(x \cdot m) =  \left\{
\begin{array}{cc}
1 & m = (0, \dots , 0) \\
0 & otherwise
\end{array}\right.
\]
\end{lemma}
In turn Lemma \ref{orthog} gives Plancherel and inversion-like identities.
\begin{proposition}\label{trans}
Let $f$ and $g$ be complex-valued functions defined on $\mathbb{Z}_q^d$.  Then,
\begin{align}
f(x) =  \sum_{m \in \mathbb{F}_q^d} \chi(x \cdot m) \widehat{f}(m)
\\
q^{-d} \sum_{x \in \mathbb{Z}_q^d} f(x) \overline{g(x)} = \sum_{m \in \mathbb{Z}_q^d} \widehat{f}(m) \overline{\widehat{g}(m)}
\end{align}
\end{proposition}

\section{Proof of Theorem \ref{main}}
The proof of Theorem \ref{main} follows a similar approach as that in \cite{CIP}.  We write
\[
\nu(t) = |\{(x,y) \in E \times E : \| x - y \| = t\}|,
\]
and we will demonstrate that $\nu(t) > 0$ for each $t \in \mathbb{Z}_q$.  To this end write
\begin{align*}
\nu(t) &= \sum_{x, y} E(x) E(y) S_t(x - y)
\\
&= \sum_{x,y,m} E(x) E(y) \widehat{S_t}(m) \chi(m \cdot (x - y))
\\
&= q^{2d} \sum_{m} \left| \widehat{E}(m) \right|^2 \widehat{S_t}(m)
\\
&= q^{-d} |E|^2 |S_t| + q^{2d} \sum_{m \neq 0} \left| \widehat{E}(m) \right|^2 \widehat{S_t}(m)
\\
&= M + R_t.
\end{align*}

We will utilize the following Lemmas.

\begin{lemma}\label{spheres}
For $d > 2$ and $t \in \mathbb{Z}_q$ for odd $q$, we have
\[
|S_t| = q^{d-1}(1 + \underline{o}(1)). 
\]
\end{lemma}

\begin{lemma}\label{decay}
Let $d > 2$, and $q = p_1^{\alpha_1} \dots p_k^{\alpha_k}$, where $q$ is odd.  Then for $m \neq 0$, we have
\[
|\widehat{S}_t(m)| \leq q^{-1} \tau(q) p_{1}^{- \frac{d-2}{2}}.
\]

\end{lemma}

Applying Lemma \ref{spheres} it is immediate that
\begin{align*}
M = q^{-1} |E|^2 (1 + \underline{o}(1)).
\end{align*}
In order to deal with the error term $R_t$, we note that
\begin{align*}
|R_t| \leq q^{2d} \max_{m \neq 0} \left| \widehat{S}_t(m) \right| \sum_{m \neq 0} \left| \widehat{E}(m) \right|^2 \leq q^{d} |E| \cdot \max_{m \neq 0} \left| \widehat{S}_t(m) \right|
\end{align*}
where the last inequality follows from adding back the zero element and applying Proposition \ref{trans}.   Applying Lemma \ref{decay} and putting the estimates for $M$ and $R_t$ together, we see that
\[
\nu(t) = q^{-1} |E|^2 (1 + \underline{o}(1)) + R_t,
\]
where
\begin{align*}
|R_t|  \leq |E| \cdot \tau(q) q^{d-1} p_1^{- \frac{d-2}{2}}
\end{align*}
and this shows that $\nu(t) > 0$ whenever
\[
|E| \geq C \tau(q) q^d p_1^{- \frac{d-2}{2}}
\]
for a sufficiently large constant $C$.  It remains to prove Lemmas \ref{spheres} and \ref{decay}.

\subsection{Gauss Sums}
Before we prove Lemmas \ref{spheres} and \ref{decay}, we will need the following well known result which we provide for completeness.
\begin{definition}[Quadratic Gauss sums] \label{Def:gs}
For positive integers $a, b, n$, we denote by $G(a,b,n)$ the following sum
\[
G(a,b,n) := \sum_{x \in \mathbb{Z}_n} \chi(ax^2 + bx).
\]
where $\chi(x) = e^{2 \pi i x /n}$.  For convenience, we denote the sum $G(a,0,n)$ by $G(a,n)$.
\end{definition}
\begin{proposition}[\cite{IK04}]\label{Prop:quadgausssum}
Let $\chi(x) = e^{2 \pi i x /n}$.  For $a \in \mathbb{Z}_n$ with $(a,n) = 1$, we have
\[
G(a,n) = 
\left\{
\begin{array}{lcl}
\varepsilon_n \left( \frac{a}{n} \right) \sqrt{n}  &   & n \equiv 1 \pmod{2}  \\
0 &   & n \equiv 2 \pmod{4}  \\
(1+i)\varepsilon_a^{-1} \left( \frac{n}{a} \right) \sqrt{n}  &   &   n \equiv 0 \pmod{4}, ~ a \equiv ~1 \pmod{2}
\end{array}
\right.
\]
where $\big(\cdot\big)$ denotes the Jacobi symbol and
\[
\varepsilon_n = 
\left\{
\begin{array}{lcc}
1    &  & n \equiv 1 \pmod{4}  \\
i     &  & n \equiv 3 \pmod{4}
\end{array}
\right.
\]
Furthermore, for general values of $a \in \mathbb{Z}_n$, we have
\[
G(a,b,n) = \left\{
\begin{array}{lc}
(a,n)G\left(\frac{a}{(a,n)}, \frac{b}{(a,n)},\frac{n}{(a,n)} \right) & (a,n) \mid b \\
0 & otherwise
\end{array}\right.
\]
\end{proposition}

\subsection{Proof of Lemma \ref{spheres}}
We first note that by the Chinese Remainder Theorem, it is enough to prove Lemma \ref{spheres} in the case that $q = p^{\ell}$ is a power of a prime.  Write
\begin{align*}
|S_t| &= \sum_{x \in \mathbb{Z}_q^d} S_t(x)
\\
&= q^{-1} \sum_{x \in \mathbb{Z}_q^d} \sum_{s \in \mathbb{Z}} \chi(s (x_1^2 + \dots + x_d^2 - t))
\\
&= q^{d-1} + q^{-1} \sum_{s \neq 0} \chi(-st)   \left( G(s, p^{\ell})\right)^d
\\
&= q^{d-1} + II_t.
\end{align*}
Put $\text{val}_p(s) = k$, if $p^k$ is the largest power of $p$ dividing $s$, in which case we let $s = p^k u$, where $u \in \mathbb{Z}_{p^{\ell - k}}^{\times}$ is uniquely determined.  Then, 
\begin{align*}
II_t  &= q^{-1} \sum_{k = 0}^{\ell - 1} \sum_{\text{val}_p(s) = k} \chi(-st) \left( G(s, p^{\ell})\right)^d
\\
&= q^{-1} \sum_{k = 0}^{\ell - 1} \sum_{u \in \mathbb{Z}_{p^{\ell - k}}^{\times} }\chi(-p^k u t) \left( G(p^k u, p^{\ell}) \right)^d
\\
&= q^{-1} \sum_{k = 0}^{\ell - 1} p^{kd} \sum_{u \in \mathbb{Z}_{p^{\ell - k}}^{\times}} \chi( - p^k ut) \left( G(u, p^{\ell - k}) \right)^d
\\
&= q^{-1} \sum_{k = 0}^{\ell - 1} p^{kd}  \left(p^{\frac{\ell - k}{2}} \right)^d \sum_{u \in \mathbb{Z}_{p^{\ell - k}}^{\times}} \chi(- p^k ut) \left( \frac{u}{p} \right)^{d (\ell -k)}.
\end{align*}
Applying the trivial bound
\[
\left| \sum_{u \in \mathbb{Z}_{p^{\ell - k}}^{\times}} \chi(- p^k ut) \left( \frac{u}{p} \right)^{d (\ell -k)} \right| \leq  p^{\ell - k},
\] we have
\begin{align*}
|II_t| &\leq q^{-1} \sum_{k=0}^{\ell - 1} p^{kd} p^{\left( \frac{\ell - k}{2} \right)d} p^{\ell - k}
\\
&= p^{\frac{\ell d}{2}} \sum_{k=0}^{\ell - 1} p^{ \frac{ k(d-2)}{2}}
\\
&\leq \ell p^{\frac{\ell d}{2}} p^{\frac{(\ell - 1)(d-2)}{2}}
\\
&= \ell q^{\left(\frac{2 \ell - 1}{2 \ell} \right)d} q^{\frac{1}{\ell} - 1}
\\
&= \ell q^{d-1} q^{\frac{1}{\ell} \left( 1 - \frac{d}{2} \right)}
\end{align*}
which shows that $II_t = \underline{o}(q^{d-1})$ whenever $d > 2$.

\subsection{Proof of Lemma \ref{decay}}
Writing $m = (m_1, \dots , m_d)$ and unraveling the definition, we see
\begin{align*}
\widehat{S}_t(m) &= q^{-d-1} \sum_{s \in \mathbb{Z}_q} \sum_{x \in \mathbb{Z}_q^d} \chi(s(x_1^2 + \dots + x_d^2 - t)) \chi(- m \cdot x)
\\
&= q^{-d-1} \sum_{s \in \mathbb{Z}_q \setminus \{0\}} \sum_{x \in \mathbb{Z}_q^d} \chi(-st) \chi(s x_1^2 - m_1x_1)\dots \chi(s x_d^2 - m_dx_d)
\\
&= q^{-d-1} \sum_{s \in \mathbb{Z}_q \setminus \{0\}} \chi(-st) \prod_{i = 1}^d G(s, - m_i, q),
\end{align*}
Our first step is to write $s = p_1^{\beta_1} \dots p_k^{\beta_k}u$ where $u \in \mathbb{Z}_{q'}^{\times}$ is uniquely determined for $q' = p_1^{\alpha_1 - \beta_1} \dots p_k^{\alpha_k - \beta_k}$
and $\beta_i \geq 0$.  Note that $\beta_i < \alpha_i$ for some $i$ since $s \neq 0$.    We will use the notation $\sum_{\beta}$ to denote the sum over all such $(\beta_1, \dots, \beta_k)$.  For $m = (m_1, \dots , m_d)$ and $\beta = (\beta_1, \dots , \beta_k)$, we define $\lambda_{m,\beta}$ to be $1$ if $p_1^{\beta_1} \dots p_k^{\beta_k} \mid m_i$ for all $i$, and zero otherwise.  When $\lambda_{m,\beta} = 1$, we put $\mu_i = \frac{m_i}{p_1^{\beta_1} \dots p_k^{\beta_k}}$.  Applying Proposition \ref{Prop:quadgausssum},
\begin{align*}
\widehat{S}_t(m) &= q^{-d-1} \sum_{\beta} \sum_{u \in \mathbb{Z}_{q'}^{\times}} \chi(-p_1^{\beta_1} \dots p_k^{\beta_k}ut) \prod_{i=1}^d G\left(p_1^{\beta_1} \dots p_k^{\beta_k} u, - m_i, q\right)
\\
&= \lambda_{m,\beta} q^{-d-1} \sum_{\beta} \sum_{u \in \mathbb{Z}_{q'}^{\times}} p_1^{\beta_1d} \dots p_k^{\beta_kd} \chi(- p_1^{\beta_1} \dots p_k^{\beta_k}ut) \prod_{i =1}^d G(u, -\mu_i, q')
\\
&= \lambda_{m,\beta} q^{-d-1} \sum_{\beta} \sum_{u \in \mathbb{Z}_{q'}^{\times}} p_1^{\beta_1d} \dots p_k^{\beta_kd} \chi(- p_1^{\beta_1} \dots p_k^{\beta_k}ut) \prod_{i =1}^d \sum_{x \in \mathbb{Z}_{q'}} \chi\left(u \left(x - \frac{\mu_i}{2u} \right)^2\right) \chi\left( \frac{- \mu_i^2} {4u} \right)
\\
&= \lambda_{m,\beta} q^{-d-1} \sum_{\beta} \sum_{u \in \mathbb{Z}_{q'}^{\times}} p_1^{\beta_1d} \dots p_k^{\beta_kd} \chi(- p_1^{\beta_1} \dots p_k^{\beta_k}ut) \prod_{i =1}^d \sum_{x \in \mathbb{Z}_{q'}} \chi\left(ux^2\right) \chi\left( \frac{- \mu_i^2} {4u} \right)
\\
&= \lambda_{m,\beta} q^{-d-1} \sum_{\beta} \sum_{u \in \mathbb{Z}_{q'}^{\times}} p_1^{\beta_1 d} \dots p_k^{\beta_k d} q'^{\frac{d}{2}} \varepsilon_{q'}^d  \chi\left(- p_1^{\beta_1} \dots p_k^{\beta_k}ut - \frac{\| \mu \|}{4 u}\right) \left( \frac{u}{q'} \right)^{d}
\\
&= \lambda_{m, \beta} q^{-d-1} \sum_{\beta} p_1^{\left(\frac{\alpha_1 + \beta_1}{2} \right)d} \dots p_k^{\left( \frac{\alpha_k + \beta_k}{2} \right)d} \varepsilon_{q'}^d \sum_{u \in \mathbb{Z}_{q'}^{\times}}\chi\left(- p_1^{\beta_1} \dots p_k^{\beta_k}ut - \frac{\| \mu \|}{4 u}\right) \left( \frac{u}{q'} \right)^{d}
\end{align*}
Applying the trivial bound to the sum in $u \in \mathbb{Z}_{q'}^{\times}$, we see that
\begin{align*}
\left|\widehat{S}_t(m)\right| &\leq q^{-d-1} \sum_{\beta} p_1^{\left(\frac{\alpha_1 + \beta_1}{2} \right)d} \dots p_k^{\left( \frac{\alpha_k + \beta_k}{2} \right)d}  p_1^{\alpha_1 - \beta_1} \dots p_k^{\alpha_k - \beta_k}
\\
&= q^{-d-1} \sum_{\beta} \prod_{i=1}^k p_i^{\left(\frac{\alpha_i + \beta_i}{2} \right)d + \alpha_i - \beta_i}.
\end{align*}
Writing $\beta_i = \alpha_i - \epsilon_i$, we have
\begin{align*}
\left|\widehat{S}_t(m)\right| &\leq q^{-d-1} \sum_{\beta} \prod_{i=1}^k p_i^{ \left( \frac{2 \alpha_i - \epsilon_i}{2} \right) d + \epsilon_i}
\\
&= q^{-d-1} \sum_{\beta} \prod_{i=1}^k p_i^{\left( \alpha_i - \frac{\epsilon_i}{2} \right)d + \epsilon_i}
\\
&= q^{-d-1} \sum_{\beta}  q^d \prod_{i=1}^k p_i^{-\frac{ d \epsilon_i - 2}{2}}
\\
&\leq q^{-1} \tau(q) p_i^{- \frac{d-2}{2}}
\end{align*}
for some $p_i$ since $\epsilon_i > 0$ for at least one value $i$ as $s \neq 0$.  The largest value obtained by the quantity occurs when $\epsilon_1 = 1$ and $\epsilon_i = 0$ for $i > 1$.  Hence,
\[
\left|\widehat{S}_t(m) \right| \leq q^{-1} \tau(q) p_1^{- \frac{d-2}{2}}.
\]

\end{document}